\documentclass[twocolumn,fleqn]{autart}    % Enable this line and disable the
\usepackage{graphicx}          % Include this line if your
\usepackage{amssymb}
\usepackage[cmex10]{amsmath}
\usepackage{accents}
\usepackage{epsfig}
\usepackage{amsmath}
\usepackage{color}
\usepackage{natbib}
\usepackage{array}
\usepackage{tikz}
\usetikzlibrary{arrows,shapes,patterns,positioning,calc}
\tikzstyle{boxstyle}=[draw=black,inner sep=7pt]
\tikzstyle{arrowstyle}=[]
\usetikzlibrary{arrows, decorations.markings}

\newlength\figureheight
\newlength\figurewidth

\allowdisplaybreaks
\newcommand{\R}{\mathbb{R}}

\newcommand{\norm}[1]{\left|\left| #1 \right|\right|}

\begin{document}
% \baselineskip=.95\normalbaselineskip

\begin{frontmatter}
\title{Bilateral Boundary Control Design for a Cascaded Diffusion-ODE System Coupled at an Arbitrary Interior Point}

\thanks[footnoteinfo]{This paper was not presented at any IFAC
meeting. Corresponding author S. Chen}

\author[UCSD]{Stephen Chen}\ead{stc007@ucsd.edu},    % Add the
\author[UdS]{Rafael Vazquez}\ead{stc007@eng.ucsd.edu},
\author[UCSD]{Miroslav Krstic}\ead{krstic@ucsd.edu}   % (ead) as shown

\address[UCSD]{Department of Mechanical and Aerospace Engineering, University of California,
 9500 Gilman Dr., La Jolla, CA 92093-0411}
\address[UdS]{Departmento de Ingenieria Aeroespacial, Escuela Superior de Ingenieros, Universidad de Sevilla, Avda. de los Descubrimientos, 41092 Seville, Spain }  % Please supply

                             % keyword list or with the
                                          % help of the Automatica
                                          % keyword wizard

\begin{abstract}
	We present a methodology for designing bilateral boundary controllers for a class of systems consisting of a coupled diffusion equation with an unstable ODE at an arbitrary interior point. A folding transformation is applied about the coupling point, transforming the system into an ODE with an input channel consisting of two coupled diffusive actuation paths. A target system with an exponentially stable trivial solution in the sense of $L^2 \times \R^n$ is proposed, and the stability property is shown via the Lyapunov method. The stabilizing control laws are formulated via tiered Volterra transformations of the second kind, establishing an equivalence relation between the stable target system and the original plant under boundary feedback. Stability properties of the plant under feedback is inferred from the equivalence relation. The well-posedness of the backstepping transformations involved are studied, and the existence of bounded Volterra kernels is shown, constituting a sufficient condition for the invertibility of the Volterra transformations.
\end{abstract}
\end{frontmatter}

\section{Preliminaries}
\label{sec:prelim}
\subsection{Introduction}
Systems modeled by parabolic partial differential equations are relevant in many engineering and social systems, with applications in many varied fields. In \citep{doi:10.1146/annurev.bioeng.8.061505.095807}, the authors model tumor angiogenesis with a nonlinear coupled parabolic system. In \citep{Hastings1978}, predator-prey Lotka-Volterra population models are formulated and studied. On the social dynamics side, opinion dynamics (modeled via the Fischer-Kolmogorov-Petrovsky-Piskunov equation) have been analyzed in the economics field via work by \citet{opiniondynamics}. With more engineering applications, \citet{vazquez2007porous} studies flows through porous media via the parabolic equation arising from Darcy's Law. Often times, there is some control objective associated with these systems, especially that of stabilization.

Also of interest are systems that involve various couplings of infinite-dimensional and finite-dimensional systems. This subject, in the context of control design, has been explored significantly, of various coupling structures of equations of varying class. In particular, cascading structures with parabolic and hyperbolic actuation paths entering linear and nonlinear ODEs -- in the parabolic case, one can think of a ``smearing'' phenomena affecting the control input. The stabilization problem of parabolic PDEs coupled with ODEs via backstepping boundary control has been studied by \citet{krstic2009compensating}. This initial result has been extended to consider various different coupling topologies, including different boundary conditions in \citep{ANTONIOSUSTO2010284}, bidirectional coupling in \citep{TANG2011540}, and sliding mode control designs in \citep{WANG201523}.

A majority of boundary backstepping designs (in 1-D) are \emph{unilateral}, meaning a single scalar controller actuates at precisely one boundary. A wide variety of results have been developed for a broad class of systems under this paradigm. However, in higher dimensions, the analogous control design would be to only actuate at some subset of the boundary (rather than on the entire boundary surface). The fully actuated high dimensional boundary control case (studied on $n$-D ball geometry by \citet{vazquez2017boundary}) motivates the study of \emph{bilateral} control design in 1-D, which, as the name suggests, involves two scalar controls two boundary points (the boundary surface of a 1-D ball). The two controllers are coupled implicitly through the equation. Intuitively, the addition of one more controller augments the controllability of the system -- an analogy to having two hands versus one when performing tasks. Some bilateral boundary control design techniques for 1-D PDEs has been studied prior in other specific contexts: for parabolic PDEs in \citep{vazquez2016bilateral} and \citep{chen2019folding}, for heterodirectional hyperbolic PDE systems in \citep{auriol2018two}, and nonlinear viscous Hamilton-Jacobi PDE in \citep{nikos2018bilateralhj}, amongst others.

The system in question in this paper involves an unstable linear ODE coupled not at a boundary, but rather, at an interior point. Previous work by \citet{zhou2017stabilization} has studied this problem in the context of unilateral control design, employing a nontraditional Fredholm transformation technique with separable kernels. This is in contrast to the work proposed in this paper, which utilizes a methodology of bilateral control design called \emph{folding}. The folding approach detailed in \citep{chen2019folding} involves using a transformation to ``fold'' the system around an interior point into a coupled parabolic PDE with a degree of freedom in choosing the folding point. In this particular case, we select the coupling point to fold about to recover a type of cascaded coupled PDE-ODE system.

The ODE coupling appearing in the interior of the PDE falls in a special class of so-called ``sandwich'' systems--systems that have a tri-layer (possibly more) of systems coupled together. Certain results exist for these systems in the unilateral sense -- for example, for ODEs ``sandwiched'' between first-order hyperbolic PDEs as in work by \citet{yu2019bilateral}. This idea exists for ODEs sandwiched by parabolic equations in the work by \citet{zhou2017stabilization}. The ``opposite'' case of a parabolic equation sandwiched by ODEs is considered by \citet{8651518}. An addition to parabolic-ODE sandwich systems is related work by \citet{koga2019single} involving the two-phase Stefan problem, a special case of an ODE sandwiched by parabolic equations whose domains evolve as a function of the ODE (a nonlinear bidirectional coupling). Finally, some results from \citet{de2016boundary} also exist for ODEs sandwiched by second-order hyperbolic PDEs in the context of the Rijke tube, a phenomena found in thermoacoustics. Related results exist in work by \citet{wang2018axial}, which involves a wave equation with ODE coupling at a moving boundary. The problem considered in this paper of the heat equation with an ODE coupled at the interior point is an example of such a sandwiched system, and may be seen to be roughly analogous to a parabolic case of the linearized Rijke tube.

The paper is structured as follows: in Section \ref{sec:model}, the model is introduced and the folding transformation is applied to recover the equivalent coupled PDE-ODE system. In Section \ref{sec:feedback}, the controller is designed via applying a two-tiered backstepping approach to recover a target system with a trivial solution possessing desirable stability properties. The stability is shown via the method of Lyapunov, and the feedback controllers (in the original coordinates) are derived. In Section \ref{sec:wellposed}, the well-posedness of the transformations from Section \ref{sec:feedback} is investigated. The existence of the transformations are shown, verifying the equivalence relation between the original plant under feedback with the chosen target system. Finally, the paper is concluded in Section \ref{sec:conclusion}.

\subsection{Notation}
The partial operator is notated using the del-notation, i.e.
\begin{align*}
	\partial_x f &:= \frac{\partial f}{\partial x}
\end{align*}

We will consider several different spaces and their Cartensian products. $\R^n$ is the standard real $n$-dimensional space. An element $v \in \R^n$ has elements notated $v_i, i \in \{1,...,n\}$. The $p$-norm denoted
\begin{align*}
	|v|_p := \left( \sum_{k=1}^{n} |v_k|^p \right)^{\frac{1}{p}}
\end{align*}
We also consider the space of square-integrable functions $L^2(I)$ over two different closed intervals $I$. For notational compactness, we label the spaces $L^2(I)$ as merely $L^2$, where the domain is implicit in the function considered. The $L^2$ space is endowed with the norm
\begin{align*}
	\norm{f}_{L^2} := \bigg(\int_I |f|_2^2 \,\,d\mu \bigg)^{\frac{1}{2}}
\end{align*}
The Cartesian product space $L^2 \times \R^n$ induces a norm $\norm{\cdot}$
\begin{align*}
	\norm{(f,v)} := \sqrt{\norm{f}_{L^2}^2 + |v|_2^2}
\end{align*}

Elements of a matrix $A$ are denoted by $a_{ij}$, in reference to the $i$-th row and $j$-th column.

\section{Model and problem formulation}
\label{sec:model}
We consider the following coupled PDE-ODE system consisting of a diffusion PDE with an unstable ODE:
\begin{align}
	\partial_t u(y,t) &= \varepsilon \partial_{y}^2 u(y,t) \label{eq:model_org_1} \\
	\dot{Z}(t) &= AZ(t) + B u(y_0,t) \label{eq:model_org_2}\\
	u(-1,t) &= \mathcal{U}_1(t) \label{eq:model_org_3}\\
 	u(1,t) &= \mathcal{U}_2(t) \label{eq:model_org_4}
\end{align}
with solutions $u: [-1,1] \times [0,\infty) \rightarrow \R, Z: [0,\infty) \rightarrow \R^n$. It is assumed that $\varepsilon > 0$ for well-posedness. The controllers operate at $x = 1$ and $x = -1$, and are denoted $\bar{\mathcal{U}}_1(t),\bar{\mathcal{U}}_2(t)$, respectively. The ODE \eqref{eq:model_org_2} is forced by the state of the heat equation at an interior point $y_0 \in (-1,1)$, which is assumed to be known \emph{a priori}. The pair $(A,B)$ is assumed to be stabilizable.

In general, a general class of reaction-advection-diffusion equations with spatially varying advection and reaction can be chosen rather than the pure heat equation, i.e. equations of the form
\begin{align*}
		\partial_t u(y,t) &= \varepsilon \partial_{y}^2 u(y,t) + b(y) \partial_y u(y,t) + \lambda(y) u(y,t)
	\end{align*}
For clarity in the paper, we merely use the pure heat equation, but the analysis is analogous to the work in \cite{chen2019folding}.

\begin{figure}[tb]
	  \centering
		\includegraphics[width=\linewidth]{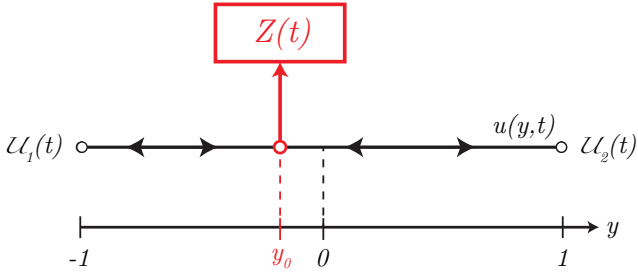}
		\caption{System schematic of heat equation coupled with interior ODE with two boundary inputs. The ODE system $Z(t)$ is located at some arbitrary interior point $y_0$.}
		\label{fig:system_schem}
\end{figure}

We perform a folding transformation about $y_0$, in which the scalar parabolic PDE system $u$ is ``folded'' into a $2 \times 2$ coupled parabolic system. We define the the folding spatial transformations as
\begin{align}
	x &= (y_0 - y)/(1 + y_0) & y \in (-1,y_0) \\
	x &= (y - y_0)/(1 - y_0) & y \in (y_0,1)
\end{align}
admitting the following states:
\begin{align}
	U(x,t) &:= \begin{pmatrix} u_1(x,t) \\ u_2(x,t) \end{pmatrix} = \begin{pmatrix} u(y_0 - (1+y_0)x,t) \\ u(y_0 + (1-y_0)x,t) \end{pmatrix} \label{eq:fold_tfm}
\end{align}
whose dynamics are governed by the following system:
\begin{align}
	\partial_t U(x,t) &= E \partial_{x}^2 U(x,t) \label{eq:model_fold_1} \\
	\dot{Z}(t) &= AZ(t) + B \Theta U(0,t) \label{eq:model_fold_2}\\
	\alpha U_x(0,t) &= -\beta U(0,t) \label{eq:model_fold_3}\\
	U(1,t) &= \mathcal{U}(t) \label{eq:model_fold_4}
\end{align}
with the parameters given by :
\begin{align}
	E &:= \text{diag}(\varepsilon_1,\varepsilon_2) \nonumber\\ &:= \text{diag}\left(\frac{\varepsilon}{(1+y_0)^2},\frac{\varepsilon}{(1-y_0)^2}\right) \\
	\alpha &:= \begin{pmatrix} 1 & a \\ 0 & 0 \end{pmatrix} \\
	\beta &:= \begin{pmatrix} 0 & 0 \\ 1 & -1 \end{pmatrix} \\
	a &:= (1+y_0)/(1-y_0) \label{eq:a_def} \\
	\Theta &= \begin{pmatrix} \theta & 1 - \theta \end{pmatrix}, \theta \in [0,1]
\end{align}
In particular, the boundary conditions \eqref{eq:model_fold_3} are curious. While they may initially appear to be Robin boundary conditions, they actually encapsulate compatibility conditions arising from imposing continuity in the solution at the folding point. Some related conditions have been considered in some previous parabolic backstepping work by \cite{tsubakino2013boundary}, albeit with differing context. This is contrasted with more typical boundary conditions which impose a single condition at a single boundary.

\begin{figure}[tb]
	  \centering
		\includegraphics[width=\linewidth]{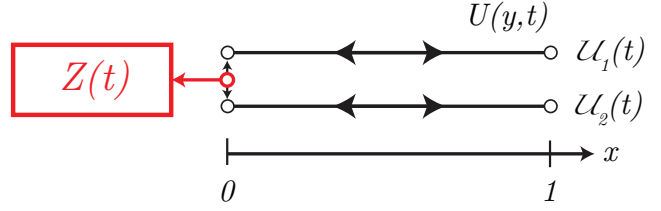}
		\caption{System schematic of folded system. The system becomes equivalent to a coupled parabolic PDE system with folding conditions imposed at the distal boundary. The folding conditions also enter the ODE as an input.}
		\label{fig:system_schem_fold}
\end{figure}

From Figure \ref{fig:system_schem_fold}, it is quite clear to see the control problem after folding becomes equivalent to stabilizing an ODE system through a coupled parabolic PDE actuation path, however, one which has the distal end ``pinned'' together. The control designs for $\mathcal{U}_{1,2}$ will be coupled.

\begin{assum}
	The ODE location $y_0$ is restricted to interval $(-1,0]$ without loss of generality. The case $y_0 \in [0,1)$ can be recovered by using a change in spatial variables $\hat{y} = -y$ and performing the same folding technique. By choosing $y_0$ in this manner, we impose an ordering $\varepsilon_1 > \varepsilon_2$.
	\label{assum:order}
\end{assum}

% \begin{figure}[tb]
	  % \centering
		% \includegraphics[width=\linewidth]{matlab/orig_sys_schem.eps}
		% \caption{System schematic of diffusion-reaction equation with two boundary inputs. The control folding point $y_0$ and the measurement location $\hat{y}_0$ can be arbitrarily chosen on the interior, independent of one another.}
		% \label{fig:system_schem}
% \end{figure}

\section{State-feedback design}
\label{sec:feedback}
The backstepping state-feedback control design is accomplished with two PDE backstepping steps. First, we will assume the existence of a stabilizing nominal control.
\begin{assum}
	There exists $\Gamma_0 \in \R^{1 \times n}$ such that the matrix $A + B \Gamma_0$ is Hurwitz.
	\label{assum:nominal_gain}
\end{assum}
Assumption \ref{assum:nominal_gain} is a direct consequence of the stabilizability of the pair $(A,B)$.

\subsection{First transformation $K$}
The first PDE backstepping transformation is a $2 \times 2$ Volterra integral transformation of the second kind:
\begin{align}
	W(x,t) &= U(x,t) - \int_0^x K(x,y) U(y,t) dy - \Gamma (x) Z(t) \label{eq:tfm1}
\end{align}
where $K \in C(\mathcal{T})$, with $\mathcal{T} := \{ (x,y) \in \R^2 | 0 \leq y \leq x \leq 1 \}$ and $\Gamma:[0,1] \rightarrow \R^{2 \times n}$. We suppose the row elements of $\Gamma$ are denoted with the index $i = 1,2$, i.e.
\begin{align}
	\Gamma(x) &:= \begin{pmatrix} \Gamma_1(x) \\ \Gamma_2(x) \end{pmatrix}
\end{align}
where $\Gamma_1(x),\Gamma_2(x) \in \R^{1 \times n}$. The associated inverse transformation is given by
\begin{align}
	U(x,t) &= W(x,t) - \int_0^x \bar{K}(x,y) W(y,t) dy - \bar{\Gamma}(x) Z(t) \label{eq:tfm1inv}
\end{align}
with $\bar{K} \in C(\mathcal{T})$ and $\bar{\Gamma}:[0,1] \rightarrow \R^{2 \times n}$. The corresponding target system for \eqref{eq:tfm1} is chosen to be
\begin{align}
	\partial_t W(x,t) &= E \partial_x^2 W(x,t) + G[K](x) W(x,t)  \label{eq:targ_int_1} \\
	\dot{Z}(t) &= (A+B\Gamma_0) Z(t) + B \Theta W(0,t) \label{eq:targ_int_2}\\
	\alpha \partial_x W(0,t) &= -\beta W(0,t) \label{eq:targ_int_3}\\
	W(1,t) &= \mathcal{V}(t) \label{eq:targ_int_4}
\end{align}
where $\mathcal{V}(t) = \begin{pmatrix} 0 & \nu_2(t) \end{pmatrix}^T$ is an auxiliary control which is designed later in the paper. The controller $\mathcal{U}(t)$ can be expressed as an operator of $\mathcal{V}(t)$ by evaluating \eqref{eq:tfm1} for $x = 1$:
\begin{align}
	\mathcal{U}(t) &:= \mathcal{V}(t) + \int_0^1 K(1,y) U(y,t) dy \label{eq:control_aux}
\end{align}
The matrix-valued operator $G[\cdot](x)$ acting on $K$ is given by
\begin{align}
	G[K](x) &= \begin{pmatrix} 0 & 0 \\ (\varepsilon_2 - \varepsilon_1) \partial_y k_{21}(x,x) & 0 \end{pmatrix} =: \begin{pmatrix} 0 & 0 \\ g[k_{21}](x) & 0 \end{pmatrix} \label{eq:G_op_def}
\end{align}
From enforcing conditions \eqref{eq:model_fold_1}-\eqref{eq:model_fold_4}, \eqref{eq:targ_int_1}-\eqref{eq:targ_int_4}, the following cascaded ODE-PDE kernel can be recovered from \eqref{eq:tfm1}:
\begin{align}
	E \partial_x^2 K(x,y) - \partial_y^2 K(x,y) E &= G[K](x) K(x,y) \label{eq:kernel_pde} \\
	E \Gamma''(x) - \Gamma(x) A + G[K](x) \Gamma(x) &= 0 \label{eq:kernel_ode}
\end{align}
subject to boundary conditions
\begin{align}
	E K(x,x) - K(x,x) E &= 0 \\
	E \partial_x K(x,x) + \partial_y K(x,x) E + E \frac{d}{dx} K(x,x) &= G[K](x) \\
	K(x,0) E \partial_x U(0) &= 0 \label{eq:K_kernel_fold1} \\
	(\Gamma(x) B \Theta - \partial_y K(x,0) E) U(0) &= 0 \label{eq:K_kernel_fold2}\\
	\Gamma(0) &= \begin{pmatrix} \Gamma_0 \\ \Gamma_0 \end{pmatrix} \label{eq:gamma_ic1} \\
	\Gamma'(0) &= \begin{pmatrix} 0 \\ 0 \end{pmatrix} \label{eq:gamma_ic2}
\end{align}
The intial condition \eqref{eq:gamma_ic1} arises from two conditions on $W(0)$: \eqref{eq:targ_int_2}, \eqref{eq:targ_int_3}. Evaluating \eqref{eq:tfm1} at $x = 0$ admits
\begin{align}
	W(0,t) &= U(0,t) - \Gamma(0) Z(t) \label{eq:w_at0}
\end{align}
From \eqref{eq:model_fold_2} and \eqref{eq:targ_int_2}, we recover a condition on $\Gamma(0)$:
\begin{align}
	\Theta \Gamma(0) &= \Gamma_0 \label{eq:gamma_0_cond1}
\end{align}
Additionally, from \eqref{eq:targ_int_3},\eqref{eq:w_at0}, we can note
\begin{align}
	\Gamma_1(0) &= \Gamma_2(0) \label{eq:gamma_0_cond2}
\end{align}
\eqref{eq:gamma_0_cond1},\eqref{eq:gamma_0_cond2} uniquely determine \eqref{eq:gamma_ic1}. The conditions \eqref{eq:gamma_ic2} are derived in an analogous manner.

A symmetry with the plant is observed with \eqref{eq:K_kernel_fold1},\eqref{eq:K_kernel_fold2} encapsulating folding conditions on $K$. By imposing \eqref{eq:model_fold_3} onto \eqref{eq:K_kernel_fold1},\eqref{eq:K_kernel_fold2}, one can recover the scalar conditions
\begin{align}
	\varepsilon_1 k_{11}(x,0) - a \varepsilon_2 k_{12}(x,0) &= 0 \\
	\varepsilon_1 k_{21}(x,0) - a \varepsilon_2 k_{22}(x,0) &= 0 \\
	\varepsilon_1 \partial_y k_{11}(x,0) + \varepsilon_2 \partial_y k_{12}(x,0) &= \Gamma_1(x) B \\
	\varepsilon_1 \partial_y k_{21}(x,0) + \varepsilon_2 \partial_y k_{22}(x,0) &= \Gamma_2(x) B
\end{align}

\subsection{Second transformation $(p,q)$}
A second transformation is designed to compensate for the term $G[K](x)$ in \eqref{eq:targ_int_1}. One can see this as the correction factor needed to compensate the interaction between the two controllers. Indeed, if one inspects the structure of the operator $G$, one may note two things. Firstly, the coupling is from the faster equation (associated with $\varepsilon_1$) to the slower equation (associaterd with $\varepsilon_2$). That is, the slower equation will have additional dynamics. Secondly, the nonzero element $g$ depends on the \emph{difference} of the diffusion coefficients. For the symmetric folding (ODE located at $x = 0$) case, the coupling does not appear.

The following transformation for designing the compensation controller $\mathcal{V}(t)$ is considered:
\begin{align}
	\Omega(x,t) &= W(x,t) - \int_0^x \begin{pmatrix} 0 & 0 \\ q(x,y) & p(x-y) \end{pmatrix} W(y,t) dy \label{eq:tfm2}
\end{align}
The corresponding inverse is given by
\begin{align}
	W(x,t) &= \Omega(x,t) - \int_0^x \begin{pmatrix} 0 & 0 \\ \bar{q}(x,y) & \bar{p}(x-y) \end{pmatrix} \Omega(y,t) dy \label{eq:tfminv2}
\end{align}
We define our target system $(\Omega,Z)$ as
\begin{align}
	\partial_t \Omega(x,t) &= E \partial_x^2 \Omega(x,t) \label{eq:targ2_1} \\
	\dot{Z}(t) &= (A+B\Gamma_0) Z(t) + B \Theta \Omega(0,t) \label{eq:targ2_2} \\
	\alpha \partial_x \Omega(0,t) &= -\beta \Omega(0,t) \label{eq:targ2_3}\\
	\Omega(1,t) &= 0 \label{eq:targ2_4}
\end{align}
The transformation \eqref{eq:tfm2}, original system model \eqref{eq:model_fold_1}-\eqref{eq:model_fold_4}, and target system \eqref{eq:targ2_1}-\eqref{eq:targ2_4} will impose a set of conditions on $p,q$ that comprise a scalar nonlocal Goursat problem:
\begin{align}
	p(x) &= a^{-1} q(x,0) \label{eq:p_orig} \\
	\varepsilon_2 \partial_x^2 q(x,y) - \varepsilon_1 \partial_y^2 q(x,y) &= (c_2 - c_1) q(x,y) \nonumber\\&\quad+  g[k_{21}](y) p(x-y) \label{eq:q_orig}
\end{align}
subject to the following boundary conditions:
\begin{align}
	\partial_y q(x,x) &= \frac{g[k_{21}](x)}{\varepsilon_2 - \varepsilon_1} \\
	q(x,x) &= 0 \label{eq:pq_bc1}\\
	\partial_y q(x,0) &= a^2 p'(x) = a \partial_x q(x,0) \label{eq:pq_bc2}
\end{align}
The kernel equations for the inverse kernels $\bar{p},\bar{q}$ are similar to those of $p,q$ respectively:
\begin{align}
	\bar{p}(x) &= a^{-1} \bar{q}(x,0) \\
	\varepsilon_2 \partial_x^2 \bar{q}(x,y) - \varepsilon_1 \partial_y^2 \bar{q}(x,y) &= (c_1 - c_2) \bar{q}(x,y) \nonumber\\&\quad- g[k_{21}](x) \bar{p}(x-y)
\end{align}
with boundary conditions
\begin{align}
	\partial_y \bar{q}(x,x) &= \frac{g[k_{21}](x)}{\varepsilon_1 - \varepsilon_2} \\
	\bar{q}(x,x) &= 0 \\
	\partial_y \bar{q}(x,0) &= -a^2 \bar{p}'(x) = -a \partial_x \bar{q}(x,0) \label{eq:pqbar_bc}
\end{align}
The PDE \eqref{eq:p_orig},\eqref{eq:q_orig} and associated boundary conditions \eqref{eq:pq_bc1},\eqref{eq:pq_bc2} are studied in previous work on folding bilateral control. The controller $\mathcal{V}(t)$ can be computed by evaluating \eqref{eq:tfm2} at $x = 1$ and using the appropriate boundary conditions:
\begin{align}
	 \mathcal{V}(t) &= \begin{pmatrix} 0 \\ \nu_2(t) \end{pmatrix}= \int_0^1 \begin{pmatrix} 0 & 0 \\ q(1,y) & p(1-y) \end{pmatrix} W(y,t) dy
\end{align}

\subsection{Stability of target system $(\Omega,Z)$}
\begin{lem}
	\label{lem:stability_targtfm}
	The trivial solution $(\Omega,Z) \equiv 0$ of the target system \eqref{eq:targ2_1}-\eqref{eq:targ2_4} is exponentially stable in the sense of the $L^2 \times \R^n$ norm. That is, there exist constants $\Pi, \mu > 0$ such that
	\begin{align}
		\norm{(\Omega(\cdot,t),Z(t))}\leq \Pi \exp \left( -\mu(t-t_0) \right) \norm{(\Omega(\cdot,t_0),Z(t_0))} \label{eq:target_bound}
	\end{align}
\end{lem}
\begin{pf}
	The proof of Lemma \ref{lem:stability_targtfm} is relatively straightforward. First consider the a Lyapunov function of the form
	\begin{align}
		V(\Omega(\cdot,t),Z) &= Z(t)^T P Z(t) + \int_0^1 \left[ \Omega(x,t)^T M \Omega(x,t) \right] dx \label{eq:lyap_func}
	\end{align}
	where $M = \textrm{diag}(a^3 m ,m), m > 0$ is an analysis parameter to be chosen later, and $P \succ 0$ is the (symmetric) solution to the Lyapunov equation
	\begin{align}
		P (A + B \Gamma_0) + (A + B \Gamma_0)^T P &= - Q \label{eq:lyap_equation}
	\end{align}
	for a chosen $Q \succ 0$. The symmetric solution $P \succ 0$ is guaranteed to exist since $A + B \Gamma_0$ is designed to be Hurwitz. We note that $V(t)$ is equivalent to the $L^2 \times \R^n$ norm:
	\begin{align}
		\Pi_1 \norm{(\Omega(\cdot,t),Z(t))}^2 \leq V(t) \leq \Pi_2 \norm{(\Omega(\cdot,t),Z(t))}^2 \label{eq:lyap_equiv}
	\end{align}
	where the coefficients $\Pi_i$ are:
	\begin{align}
		\Pi_1 &= \min \{ \lambda_{\min}(P), a^3 m\} \\
		\Pi_2 &= \max \{ \lambda_{\max}(P), m \}
	\end{align}
	Differentiating \eqref{eq:lyap_func} in time, one finds
	\begin{align}
		\dot{V}(t) &\leq (\Omega(0,t)^T \Theta^T B^T + Z(t)^T (A + B \Gamma_0)^T)) P Z(t) \nonumber\\
		&\quad+ Z(t)^T P ((A+B\Gamma_0) Z(t) + B \Theta \Omega(0,t)) \nonumber\\
		&\quad+ \int_0^1 \left[ \partial_x^2 \Omega(x,t)^T E M \Omega(x,t) \right.\nonumber\\
		&\qquad\qquad\left.+ \Omega(x,t)^T M E \partial_x^2 \Omega(x,t) \right] dx
	\end{align}
	Using integration by parts and \eqref{eq:lyap_equation} will admit
	\begin{align}
		\dot{V}(t) &\leq - Z(t)^T Q Z(t) + 2 Z(t)^T P B \Theta \Omega(0,t) \nonumber\\
		&\quad- \int_0^1 \left[ 2\partial_x \Omega(x,t)^T E M \partial_x \Omega(x,t) \right] dx \nonumber\\
		&\leq -\lambda_{\textrm{min}}(Q)|Z(t)|_2^2 + 2 Z(t)^T P B \Theta \Omega(0,t) \nonumber\\
		&\quad- 2 a^3 m \varepsilon_2 \norm{\partial_x \Omega(\cdot,t)}_{L^2}^2
	\end{align}
	Applying Young's inequality,
	\begin{align}
		\dot{V}(t) &\leq - \mu_1 |Z(t)|_2^2 - 4 \mu_2 \norm{\partial_x \Omega(\cdot,t)}_{L^2}^2 \label{eq:lyap_ineq_1}
	\end{align}
	where
	\begin{align}
		\mu_1 &= \lambda_{\textrm{min}}(Q) - \delta |P|_{2,i} |B|_{2,i} \\
		\mu_2 &= \frac{1}{4} \left( 2 a^3 \varepsilon_2 m - \frac{1}{\delta} |P|_{2,i} |B|_{2,i} \right)
	\end{align}
	The analysis parameters $\delta, m_1, m_2$ must be chosen such that $\mu_{1,2} > 0$. This is easily achievable by choosing
	\begin{align}
		\delta &< \frac{\lambda_{\min}(Q)}{\lambda_{\max}(P) |B|_{2,i}} \\
		m &> \frac{{\lambda_{\max}(P) |B|_{2,i}}}{2 \delta a^3 \varepsilon_2}
	\end{align}
	Applying Young's inequality and \eqref{eq:lyap_equiv} to \eqref{eq:lyap_ineq_1}, one finds
	\begin{align}
		\dot{V}(t) &\leq - \frac{\min\{\mu_1,\mu_2\}}{\Pi_1} V(t)
	\end{align}
	which via the comparison principle admits the bound
	\begin{align}
		V(t) &\leq \exp (- 2\mu (t - t_0)) V(t_0)
	\end{align}
	where
	\begin{align}
		\mu &:= \frac{\min\{\mu_1,\mu_2\}}{2\Pi_1}
	\end{align}
	Applying the equivalence \eqref{eq:lyap_equiv} once more recovers the bound \eqref{eq:target_bound} with $\Pi = \sqrt{\Pi_2/\Pi_1}$. This completes the proof.
\end{pf}

\subsection{Main result: closed-loop stability}
\begin{thm}
The trivial solution of the system \eqref{eq:model_org_1}-\eqref{eq:model_org_4} is exponentially stable in the sense of the $L^2 \times \R^n$ norm under the pair of state feedback control laws $\mathcal{U}_1,\mathcal{U}_2$:
\begin{align}
	\begin{pmatrix} \mathcal{U}_1(t) \\ \mathcal{U}_2(t) \end{pmatrix} &= \int_{-1}^1 \begin{pmatrix} F_1(y) \\ F_2(y) \end{pmatrix} u(y,t) dy + F_3 Z(t) \label{eq:state_feedback_controls}
\end{align}
with feedback gains $F_1,F_2$ defined as
\begin{align}
	F_1(y) &= \begin{cases} (1+y_0)^{-1} k_{11} \left(1,\frac{y_0-y}{1+y_0}\right) & y \leq y_0 \\
													(1-y_0)^{-1} k_{12} \left(1,\frac{y-y_0}{1-y_0}\right) & y > y_0 \end{cases}\label{eq:gain_f1}\\
	F_2(y) &= \begin{cases} (1+y_0)^{-1} h_{1}\left(\frac{y_0-y}{1+y_0}\right) & y \leq y_0 \\
													(1-y_0)^{-1} h_{2}\left(\frac{y-y_0}{1-y_0}\right) & y > y_0 \end{cases}\label{eq:gain_f2}\\
	h_1(y) &= k_{21}(1,y) + q(1,y) - \int_y^1 \bigl[p(1-z) k_{21}(z,y) \nonumber\\&\qquad + q(1,z)k_{11}(z,y) \bigr] dz \\
	h_2(y) &= k_{22}(1,y) + p(1-y) - \int_y^1 \bigl[p(1-z) k_{22}(z,y) \nonumber\\&\qquad+ q(1,z)k_{12}(z,y) \bigr] dz \\
	F_3 &= \begin{pmatrix} \Gamma_1(1) \\ \Gamma_2(1) - \int_0^1 [q(1,y) \Gamma_1(y) + p(1-y) \Gamma_2(y)] dy \end{pmatrix}
\end{align}
where $k_{ij},p,q \in C(\mathcal{T})$ are solutions to the kernel PDE equations \eqref{eq:kernel_pde},\eqref{eq:p_orig},\eqref{eq:q_orig} respectively (with associated boundary conditions), and $\Gamma_{1,2} \in C([0,1])$ are solutions to the kernel ODE equations\eqref{eq:kernel_ode}. That is, under the feedback controllers \eqref{eq:state_feedback_controls}, there exists a constant $\bar{\Pi}$ such that
	\begin{align}
		\norm{(u(\cdot,t),Z(t))} \leq \bar{\Pi} \exp \left( -\mu(t-t_0) \right) \norm{(u(\cdot,t_0),Z(t_0))} \label{eq:state_feedback_bound}
	\end{align}
	\label{thm:statefdback}
\end{thm}
The proof of Theorem \ref{thm:statefdback} is not given but is analogous to the proofs found in the work by \cite{smyshlyaev2005spacetime}. The proof involves utilizing the invertability of the transformations \eqref{eq:fold_tfm},\eqref{eq:tfm1},\eqref{eq:tfm2} that arise either trivially (folding), or from the boundedness of the kernels (studied in Section \ref{sec:wellposed}). The forward and inverse transforms give estimates on the equivalence relation between the target system \eqref{eq:targ2_1}-\eqref{eq:targ2_4} and the original system \eqref{eq:model_org_1}-\eqref{eq:model_org_4}, which are applied to \eqref{eq:state_feedback_bound}.

\section{Well-posedness of $K,\Gamma$ kernel system}
\label{sec:wellposed}
The PDE gain kernel $K$ and ODE kernel $\Gamma$ must be shown to be well-posed. The following lemmas establish these results.

\subsection{$\Gamma$ kernel}
The ODE system \eqref{eq:kernel_ode} is written into into two separate $n$-th order ODEs:
\begin{align}
	\Gamma_1''(x) &= \varepsilon_1^{-1} \Gamma_1(x) A \\
	\Gamma_2''(x) &= \varepsilon_2^{-1} \Gamma_2(x) A + \varepsilon_2^{-1} g[k_{21}](x) \Gamma_1(x)
\end{align}
where the initial conditions can be found from \eqref{eq:gamma_ic1},\eqref{eq:gamma_ic2} to be
\begin{align}
	\Gamma_1(0) &= \Gamma_0 \\
	\Gamma_1'(0) &= 0 \\
	\Gamma_2(0) &= \Gamma_0 \\
	\Gamma_2'(0) &= 0
\end{align}
From the variation of constants formula, it is easy to see that the solutions for $\Gamma_1$ can be expressed via
\begin{align}
	\Gamma_1(x) &= \begin{pmatrix} \Gamma_0 & 0 \end{pmatrix} \exp\left(\begin{pmatrix} 0 & I \\ \varepsilon_1^{-1} A & 0 \end{pmatrix} x\right)\begin{pmatrix}I \\ 0 \end{pmatrix} \label{eq:gamma1_soln}
\end{align}
while $\Gamma_2$ is
\begin{align}
	\Gamma_2(x) &= \bigg[\begin{pmatrix} \Gamma_0 & 0 \end{pmatrix} \exp\left(\begin{pmatrix} 0 & I \\ \varepsilon_2^{-1} A & 0 \end{pmatrix} x\right) \nonumber\\
	&\qquad+ \int_0^x \begin{pmatrix} 0 & \varepsilon_2^{-1} g[k_{21}](\xi) \Gamma_1(\xi) \end{pmatrix} \nonumber\\
	&\qquad\qquad\times\exp\left(\begin{pmatrix} 0 & I \\ \varepsilon_2^{-1} A & 0 \end{pmatrix} (x - \xi)\right) d\xi \bigg] \begin{pmatrix}I \\ 0 \end{pmatrix} \label{eq:gamma2_soln}
\end{align}
Noting that $\Gamma_2(x)$ depends on $g[k_{21}](x)$, an element of $K$, we define the operator formulation $(\check{\Gamma}_2 \circ g)[k_{21}] := \Gamma_2$. We will additionally define the operator $\Phi: C(\mathcal{T};\R) \rightarrow C([0,1];\R^{2 \times n})$ as
\begin{align}
	\Phi[f](x) &:= \begin{pmatrix} \Gamma_1(x) \\ (\check{\Gamma}_2 \circ g)[f](x) \end{pmatrix}
\end{align}
where $\Phi[f]$ maps from the scalar $C(\mathcal{T};\R)$ to the multidimensional $C([0,1];\R^{2 \times n})$ function space. From this definition, it naturally follows that $\Phi[k_{21}] = \Gamma$. This operator representation will be used in the well-posedness analysis for the PDE kernel $K$.

Of interest are bounds on $\Gamma_1(x), \Gamma_2(x)$. The following bound on $\Gamma_1(x)$ is trivial to find:
\begin{align}
	|\Gamma_1(x)| &\leq |\Gamma_0|_2 S(x)
\end{align}
where
\begin{align}
	S(x) &= \sigma_{\max}\left(\exp\left( \begin{pmatrix} 0 & I \\ \varepsilon_2^{-1} A & 0 \end{pmatrix}  x \right)\right)
\end{align}
and $\sigma_{\max}(X)$ denotes the largest singular value of the matrix $X$ (the induced 2-norm). It is important to note that the largest singular value of the matrix exponential is bounded on the compact set $[0,1]$, i.e. $|S(x)| < \infty, \forall x \in [0,1]$.
With a bound on $\Gamma_1$, the following bound on $(\check{\Gamma}_2 \circ g)[k_{21}](x)$ can be found, which will be used in the proof of well-posedness for $K$.
\begin{align}
	|(\check{\Gamma}_2 \circ g)[k_{21}](x)| &\leq |\Gamma_0|_2 S(x) \nonumber\\
	&\quad+ \int_0^x \frac{\varepsilon_1 - \varepsilon_2}{\varepsilon_2} \norm{\Gamma_1}_{L^\infty} S(x-\xi) \nonumber\\
	&\qquad\qquad\qquad \times |\partial_y k_{21}(\xi,\xi)| d\xi \label{eq:gam2_bound}
\end{align}

\subsection{$K$ kernel}
The $K$-kernel must be approached in two sets of equations: $(k_{11},k_{12})$ and $(k_{21},k_{22})$. The reason for this is that the operator $G[K]$ introduces coupling between the two sets of kernels. We first apply a transformation to gain kernel \eqref{eq:kernel_pde}
\begin{align}
		\check{K}(x,y) &= \sqrt{E} \partial_x K(x,y) + \partial_y K(x,y) \sqrt{E} \label{eq:kernel_tfm1}
\end{align}
which transforms the kernel PDE into a $2 \times 2$ system of coupled first-order hyperbolic PDEs.

\subsubsection{$(k_{11},k_{12})$-system}
The transform \eqref{eq:kernel_tfm1} will admit the following coupled $2 \times 2$ system
\begin{align}
	\sqrt{\varepsilon_1} \partial_x k_{11}(x,y) + \sqrt{\varepsilon_1} \partial_y k_{11}(x,y) &= \check{k}_{11}(x,y) \label{eq:k_11_kernel}\\
	\sqrt{\varepsilon_1} \partial_x k_{12}(x,y) + \sqrt{\varepsilon_2} \partial_y k_{12}(x,y) &= \check{k}_{12}(x,y) \\
	\sqrt{\varepsilon_1} \partial_x \check{k}_{11}(x,y) - \sqrt{\varepsilon_1} \partial_y \check{k}_{11}(x,y) &= 0 \label{eq:kc_11_kernel}\\
	\sqrt{\varepsilon_1} \partial_x \check{k}_{12}(x,y) - \sqrt{\varepsilon_2} \partial_y \check{k}_{12}(x,y) &= 0\label{eq:kc_12_kernel}
\end{align}
with boundary conditions
\begin{align}
	k_{11}(x,0) &= \frac{a\varepsilon_2}{\varepsilon_1(a \varepsilon_2 + \sqrt{\varepsilon_1 \varepsilon_2})} \nonumber\\&\quad\times \int_0^x [\sqrt{\varepsilon_1} \check{k}_{11}(y,0) + \sqrt{\varepsilon_2} \check{k}_{12}(y,0) \nonumber\\ &\qquad\qquad\qquad\qquad\qquad\qquad- \Gamma_1(y) B ] dy \\
	k_{12}(x,0) &= \frac{1}{a \varepsilon_2 + \sqrt{\varepsilon_1 \varepsilon_2}} \nonumber\\&\quad\times \int_0^x [\sqrt{\varepsilon_1} \check{k}_{11}(y,0) + \sqrt{\varepsilon_2} \check{k}_{12}(y,0) \nonumber\\ &\qquad\qquad\qquad\qquad\qquad\qquad - \Gamma_1(y) B ] dy  \\
	k_{12}(x,x) &= 0 \label{eq:k12_bc1}\\
	\check{k}_{11}(x,x) &= 0 \label{eq:kc11_bc1}\\
	\check{k}_{12}(x,x) &= 0 \label{eq:kc12_bc1}
\end{align}

\begin{lem}
	The system of first-order hyperbolic PDEs \eqref{eq:k_11_kernel}-\eqref{eq:kc_12_kernel} and associated boundary conditions admit a unique set of $k_{11},k_{12} \in C(\mathcal{T})$ solutions.
	\label{lem:wellposedk1}
\end{lem}
\begin{pf}
	$(k_{11},k_{12})$ can actually be solved explicitly via the method of characteristics. First, note that \eqref{eq:kc_11_kernel},\eqref{eq:kc_12_kernel},\eqref{eq:kc11_bc1},\eqref{eq:kc12_bc1} imply that $(\check{k}_{11},\check{k}_{12}) = 0$, the solution can be simplified significantly:
	\begin{align}
		k_{11}(x,y) &= \frac{a\varepsilon_2}{\varepsilon_1(a \varepsilon_2 + \sqrt{\varepsilon_1 \varepsilon_2})} \int_0^{x-y} \Gamma_1(z) B dz \\
		k_{12}(x,y) &= \begin{cases} k_{12,l}(x,y) & \sqrt{\varepsilon_2} x \geq \sqrt{\varepsilon_1}y \\ 0 & \textrm{otherwise} \end{cases} \\
		k_{12,l}(x,y) &= \frac{1}{a \varepsilon_2 + \sqrt{\varepsilon_1 \varepsilon_2}} \int_0^{x - \sqrt{\frac{\varepsilon_1}{\varepsilon_2}}y} \Gamma_1(z) B dz
	\end{align}
\end{pf}
It is not difficult to see that $k_{11} \in C^\infty(\mathcal{T}) \subset C(\mathcal{T})$ from the matrix exponential properties, while $k_{12} \in C(\mathcal{T})$.

\subsubsection{$(k_{21},k_{22})$-system}
The $(k_{21},k_{22})$ system must be treated in a differing manner than the $(k_{11},k_{12})$ system, due to the existence of a nonlocal trace term
To account for the different nature of these characteristics, we perform one more transformation on the kernels for $k_{2i}$:
\begin{align}
	\hat{k}_{2i}(x,y) &= \sqrt{\varepsilon_2} \partial_x k_{2i}(x,y) - \sqrt{\varepsilon_i} \partial_y k_{2i}(x,y) \label{eq:kh2_tfm}
\end{align}
where $i \in \{1,2\}$. We then turn our attention to the gain kernel system $(\hat{k}_{21},\check{k}_{21},\hat{k}_{22},\check{k}_{22})$.

The component system of kernel PDEs for $(\hat{k}_{21},\check{k}_{21},\hat{k}_{22},\check{k}_{22})$ is
\begin{align}
	\sqrt{\varepsilon_2} \partial_x \hat{k}_{21}(x,y) + \sqrt{\varepsilon_1} \partial_y \hat{k}_{21}(x,y) &= - g[k_{21}](x) k_{11}(x,y) \label{eq:k_21_kernel}\\
	\sqrt{\varepsilon_2} \partial_x \hat{k}_{22}(x,y) + \sqrt{\varepsilon_2} \partial_y \hat{k}_{22}(x,y) &= - g[k_{21}](x) k_{12}(x,y) \\
	\sqrt{\varepsilon_2} \partial_x \check{k}_{21}(x,y) - \sqrt{\varepsilon_1} \partial_y \check{k}_{21}(x,y) &= - g[k_{21}](x) k_{11}(x,y) \\
	\sqrt{\varepsilon_2} \partial_x \check{k}_{22}(x,y) - \sqrt{\varepsilon_2} \partial_y \check{k}_{22}(x,y) &= - g[k_{21}](x) k_{12}(x,y) \label{eq:kc_22_kernel}
\end{align}
subject to the following boundary conditions:
\begin{align}
	\hat{k}_{21}(x,0) &= -\frac{1-a^2}{1+a^2} \check{k}_{21}(x,0) +  \frac{2a^3}{1+a^2} \check{k}_{22}(x,0) \nonumber\\&\qquad\qquad\qquad- \frac{2a^2}{\sqrt{\varepsilon_1}(1+a^2)} (\check{\Gamma}_2 \circ g)[k_{21}](x) B \label{eq:k21_reflection1} \\
	\hat{k}_{22}(x,0) &= \frac{2}{a(1+a^2)} \check{k}_{21}(x,0) + \frac{1-a^2}{1+a^2} \check{k}_{22}(x,0) \nonumber\\&\qquad\qquad\qquad- \frac{2}{\sqrt{\varepsilon_2}(1+a^2)} (\check{\Gamma}_2 \circ g)[k_{21}](x) B  \label{eq:k22_reflection1}\\
	\check{k}_{21}(x,x) &= - \frac{\sqrt{\varepsilon_1} - \sqrt{\varepsilon_2}}{\sqrt{\varepsilon_1} + \sqrt{\varepsilon_2}} \hat{k}_{21}(x,x) \label{eq:k21_reflection2}\\
	\check{k}_{22}(x,x) &= 0 \label{eq:k21_bc_end}
\end{align}
where the inverse transformations are given to be
\begin{align}
	k_{21}(x,y) &= \frac{1}{2\sqrt{\varepsilon_2}} \int_y^x \check{k}_{21}(z,y) + \hat{k}_{21}(z,y) dz \label{eq:k21_invtfm}\\
	k_{22}(x,y) &= \frac{1}{2\sqrt{\varepsilon_2}} \int_y^x \check{k}_{22}(z,y) + \hat{k}_{22}(z,y) dz \label{eq:k22_invtfm}
\end{align}
and the function $g[k_{21}](x)$ can be expressed in terms of $\hat{k}_{21},\check{k}_{21}$:
\begin{align}
	g[k_{21}](x) &= \frac{(\varepsilon_2 - \varepsilon_1)}{2 \sqrt{\varepsilon_1}} (\check{k}_{21}(x,x) - \hat{k}_{21}(x,x)) \nonumber\\
	&= (\sqrt{\varepsilon_1} - \sqrt{\varepsilon_2}) \hat{k}_{21}(x,x)
\end{align}

\begin{lem}
	The system of first-order hyperbolic PDE \eqref{eq:k_21_kernel}-\eqref{eq:kc_22_kernel} and associated boundary conditions admit a unique set of $\hat{k}_{21},\check{k}_{21},\hat{k}_{22},\check{k}_{22} \in C(\mathcal{T})$ solutions.
	\label{lem:wellposedk2}
\end{lem}
\begin{pf}
	We recognize that $(\hat{k}_{21},\check{k}_{21},\hat{k}_{22},\check{k}_{22})$ are similar in structure to the (previous result), albeit with an additional non-local recirculation term appearing in the boundaries \eqref{eq:k21_reflection1},\eqref{eq:k22_reflection1}. The non-local term in the boundary does not change the method of the proof by too much, however, additional care must be given to incorporate the behavior.

	\begin{figure}[tb]
			\centering
			\includegraphics[width=\linewidth]{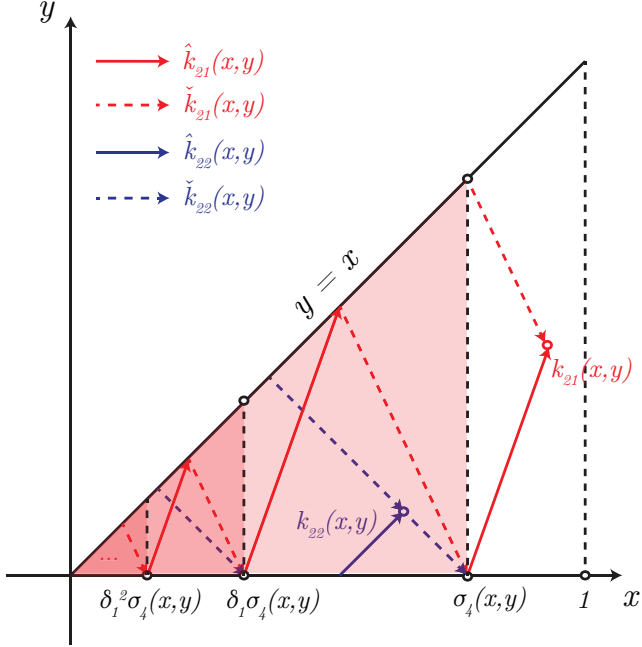}
			\caption{Characteristics of $\hat{k}_{21},\hat{k}_{22},\check{k}_{21},\check{k}_{22}$ featuring an infinite number of reflection boundary conditions. To solve for a given point, the solution must be known on a triangle of smaller volume, motivating a recursive procedure for solving the gain kernel.}
			\label{fig:k21_ode}
	\end{figure}

	The solutions for the $\hat{k}_{21}$ can be recovered via a direct application of the method of characteristics:
	\begin{align}
		\hat{k}_{21}(x,y) &= \hat{k}_{21}\left( \sigma_4(x,y), 0\right) + \hat{I}_{21}[\hat{k}_{21}](x,y) \label{eq:kh_21_inteq} \\
		\hat{k}_{22}(x,y) &= \hat{k}_{22}\left( x - y, 0\right) + \hat{I}_{22}[\hat{k}_{21}](x,y) \\
		\check{k}_{21}(x,y) &= \check{k}_{21}\left( \sigma_5(x,y), \sigma_5(x,y)\right) + \check{I}_{21}[\hat{k}_{21}](x,y) \\
		\check{k}_{22}(x,y) &= \check{k}_{22}\left( \frac{x+y}{2}, \frac{x+y}{2}\right) + \check{I}_{22}[\hat{k}_{21}](x,y) \label{eq:kc_21_inteq}
	\end{align}
	\begin{align}
		\sigma_4(x,y) &:= x-\frac{\sqrt{\varepsilon_2}}{\sqrt{\varepsilon_1}}y \\
		\sigma_5(x,y) &:= \frac{\sqrt{\varepsilon_1}x + \sqrt{\varepsilon_2}y}{\sqrt{\varepsilon_1} + \sqrt{\varepsilon_2}}
	\end{align}
	and the integral operators $\hat{I}_{21},\hat{I}_{22},\check{I}_{21},\check{I}_{22}$ are defined
	\begin{align}
		&\hat{I}_{21}[\hat{k}_{21}](x,y) \nonumber\\
		&\quad:= \int_0^{\frac{y}{\sqrt{\varepsilon_1}}} \bigg[ - \frac{\varepsilon_2 - \varepsilon_1}{\sqrt{\varepsilon_1} + \sqrt{\varepsilon_2}} k_{11}(\sqrt{\varepsilon_2} z + \sigma_4(x,y),\sqrt{\varepsilon_1}z) \nonumber\\
		&\qquad\qquad \times\hat{k}_{21}(\sqrt{\varepsilon_2} z + \sigma_4(x,y), \sqrt{\varepsilon_2} + \sigma_4(x,y)) \bigg] dz \\
		&\hat{I}_{22}[\hat{k}_{21}](x,y) \nonumber\\
		&\quad:= \int_0^{\frac{y}{\sqrt{\varepsilon_2}}} \bigg[ - \frac{\varepsilon_2 - \varepsilon_1}{\sqrt{\varepsilon_1} + \sqrt{\varepsilon_2}} k_{12}(\sqrt{\varepsilon_2} z + x-y,\sqrt{\varepsilon_2}z) \nonumber\\
		&\qquad\qquad\qquad \times\hat{k}_{21}(\sqrt{\varepsilon_2} z + x-y,\sqrt{\varepsilon_2} + x-y) \bigg] dz \\
		&\check{I}_{21}[\hat{k}_{21}](x,y) \nonumber\\
		&\quad:= \int_0^{\frac{x-y}{\sqrt{\varepsilon_1} + \sqrt{\varepsilon_2}}} \bigg[ - \frac{\varepsilon_2 - \varepsilon_1}{\sqrt{\varepsilon_1} + \sqrt{\varepsilon_2}} k_{11}(\sqrt{\varepsilon_2} z + \sigma_5(x,y)\nonumber\\
		&\qquad\qquad\qquad\qquad\qquad\qquad,-\sqrt{\varepsilon_1}z + \sigma_5(x,y)) \nonumber\\
		&\qquad\qquad\qquad\qquad \times\hat{k}_{21}(\sqrt{\varepsilon_2} z + \sigma_5(x,y), \nonumber\\
		&\qquad\qquad\qquad\qquad\qquad\qquad \sqrt{\varepsilon_2} + \sigma_5(x,y)) \bigg] dz \\
		&\check{I}_{22}[\hat{k}_{21}](x,y) \nonumber\\
		&\quad:= \int_0^{\frac{x-y}{2\sqrt{\varepsilon_2}}} \bigg[ - \frac{\varepsilon_2 - \varepsilon_1}{\sqrt{\varepsilon_1} + \sqrt{\varepsilon_2}} \nonumber\\
		&\qquad\qquad \times k_{12}\left(\sqrt{\varepsilon_2} z + \frac{x+y}{2},-\sqrt{\varepsilon_2}z + \frac{x+y}{2}\right) \nonumber\\
		&\qquad\qquad \times\hat{k}_{21}\left(\sqrt{\varepsilon_2} z + \frac{x+y}{2}, \sqrt{\varepsilon_2} + \frac{x+y}{2}\right) \bigg] dz
	\end{align}

	From enforcing \eqref{eq:k21_reflection1}-\eqref{eq:k21_bc_end} on \eqref{eq:kh_21_inteq}-\eqref{eq:kc_21_inteq} recursively, one can eventually arrive at an integral equation system representation for $(\hat{k}_{21},\hat{k}_{22},\check{k}_{21},\check{k}_{22})$ involving infinite sums of Volterra-type integral operators. The infinite sums appear due to the reflection boundary conditions \eqref{eq:k21_reflection1},\eqref{eq:k21_reflection2} observed in the system.
	\begin{align}
		&\hat{k}_{21}(x,y) = \nonumber\\
		&\quad \lim_{n \rightarrow \infty} \bigg[ - \delta_1^n\delta_2^{n+1} \check{k}_{21} \left( \delta_1^n \sigma_4(x,y),0\right) \nonumber\\
		&\qquad\qquad+ \delta_1^n \delta_2^n \frac{2 a^3}{1 + a^2} \check{k}_{22}(\delta_1^n \sigma_4(x,y),0)\bigg] \nonumber\\
		&\quad+ \sum_{n=0}^{\infty} \bigg[\delta_1^n \delta_2^n \hat{I}_{21}[\hat{k}_{21}](\delta_3^n \sigma_4(x,y), \delta_3^n \sigma_4(x,y)) \nonumber\\
		&\qquad- \delta_1^n \delta_2^{n+1} \check{I}_{21}[\hat{k}_{21}](\delta_1^n \sigma_4(x,y),0) \nonumber\\
		&\qquad+ \delta_1^n \delta_2^n \frac{2 a^3}{1 + a^2} \check{I}_{22}[\hat{k}_{21}](\delta_1^n \sigma_4(x,y),0) \nonumber\\
		&\qquad- \delta_1^n \delta_2^n \frac{2a^2}{\sqrt{\varepsilon_1}(1+a^2)} (\check{\Gamma}_2 \circ g)[k_{21}](\delta_1^n \sigma_4(x,y)) \bigg] \nonumber\\
		&\quad+ \hat{I}_{21}[\hat{k}_{21}](x,y) \label{eq:kh_21_eq}\\
		&\hat{k}_{22}(x,y) = \nonumber\\
		&\quad \frac{2}{a(1+a^2)} \lim_{n \rightarrow \infty} \bigg[ \delta_1^n \delta_2^n \check{k}_{21}(\delta_1^n (x-y),0) \bigg] \nonumber\\
		&\quad+ \frac{2}{a(1+a^2)} \sum_{n=1}^{\infty} \bigg[ (-1)^n \delta_1^n \delta_2^{n-1} \nonumber\\
		&\qquad\qquad\qquad\times\hat{I}_{21}[\hat{k}_{21}](\delta_3^n (x-y),\delta_3^n (x-y)) \nonumber\\
		&\qquad + \delta_1^{n-1} \delta_2^{n-1} \check{I}_{21}[\hat{k}_{21}](\delta_1^n (x-y),0) \nonumber\\
		&\qquad + (-1)^n \delta_1^n \delta_2^{n-1} \frac{2a}{1+a^3} \check{I}_{22}[\hat{k}_{21}](\delta_3^n (x-y),0) \nonumber\\
		&\qquad + \delta_1^{n} \delta_2^{n-1} \frac{2a^2}{\sqrt{\varepsilon_1} (1 + a^2)} (\check{\Gamma}_2 \circ g)[k_{21}](\delta_1^n (x-y)) \bigg] \nonumber\\
		&\quad+ \delta_2 \check{I}_{22}[\hat{k}_{21}](x-y,0) + \hat{I}_{22}[\hat{k}_{21}](x,y) \nonumber\\
		&\quad- \frac{2}{\sqrt{\varepsilon_2}(1+a^2)} (\check{\Gamma}_2 \circ g)[k_{21}](x-y) B \label{eq:kh_22_eq}\\
		&\check{k}_{21}(x,y) = \nonumber\\
		&\quad \lim_{n \rightarrow \infty} \bigg[\delta_1^n \delta_2^n \check{k}_{21}(\delta_1^n \sigma_5(x,y), \delta_1^n \sigma_5(x,y)) \bigg] \nonumber\\
		&\quad + \sum_{n=0}^{\infty} \bigg[ - \delta_1^{n+1} \delta_2^n \hat{I}_{21}[\hat{k}_{21}](\delta_1^n \sigma_5(x,y), \delta_1^n \sigma_5(x,y)) \nonumber\\
		&\qquad - \frac{2a^3}{1+a^2} \delta_1^{n+1} \delta_2^n \check{I}_{22}[\hat{k}_{21}]\left(\frac{\delta_1}{\delta_3}\delta_1^{n} \sigma_5(x,y),0\right) \nonumber\\
		&\qquad - \delta_1^{n+1} \delta_2^n \frac{2a^2}{\sqrt{\varepsilon_1}(1+a^2)} (\check{\Gamma}_2 \circ g)[k_{21}]\left(\frac{\delta_1}{\delta_3}\delta_1^{n} \sigma_5(x,y)\right) \nonumber\\
		&\qquad +\delta_1^{n+1} \delta_2^{n+1} \check{I}_{21}[\hat{k}_{21}]\left(\frac{\delta_1}{\delta_3}\delta_1^{n} \sigma_5(x,y),0\right) \bigg] \nonumber\\
		&\quad + \check{I}_{21}[\hat{k}_{21}](x,y) \label{eq:kc_21_eq}\\
		&\check{k}_{22}(x,y) = \check{I}_{22}[\hat{k}_{21}](x,y) \label{eq:kc_22_eq}
	\end{align}
	with
	\begin{align}
		\delta_1 &= \frac{\sqrt{\varepsilon_1} - \sqrt{\varepsilon_2}}{\sqrt{\varepsilon_1}+\sqrt{\varepsilon_2}} \\
		\delta_2 &= \frac{1-a^2}{1+a^2} \\
		\delta_3 &= \frac{\sqrt{\varepsilon_1}}{\sqrt{\varepsilon_1} + \sqrt{\varepsilon_2}}
	\end{align}
	Since $a < 1, \varepsilon_1 > \varepsilon_2$ as per Assumption \ref{assum:order}, the coefficients $\delta_{1,2,3} \in (0,1)$. It is unclear initially whether the limit and infinite sum terms even converge, however, as one may notice in Figure \ref{fig:k21_ode}, the contracting volume of integration and the reflection coefficients (appearing implicitly in the $\delta_i$ coefficients) will guarantee the convergence.

	The argument is similar to that of (previous work), but the additional non-local term $(\check{\Gamma_2} \circ g)[k_{21}]$ must be accounted for. As we have shown in \eqref{eq:gam2_bound}, the operator $\check{\Gamma}_2 \circ g)$ is a bounded Volterra operator on the transformed variable $\hat{k}_{21}$. It is easy to see that if a solution $\hat{k}_{21}$ exists, then the solutions $\check{k}_{21},\hat{k}_{22},\check{k}_{22}$ follow from direct evaluation. Thus, in this proof, we will merely show that $\hat{k}_{21}$ exists.

	We establish an iteration $\hat{k}_{21,i}$ as
	\begin{align}
		&\hat{k}_{21,i}(x,y) = \nonumber\\
		&\quad \sum_{n=0}^{\infty} \bigg[\delta_1^n \delta_2^n \hat{I}_{21}[\hat{k}_{21,i}](\delta_3^n \sigma_4(x,y), \delta_3^n \sigma_4(x,y)) \nonumber\\
		&\qquad- \delta_1^n \delta_2^{n+1} \check{I}_{21}[\hat{k}_{21,i}](\delta_1^n \sigma_4(x,y),0) \nonumber\\
		&\qquad+ \delta_1^n \delta_2^n \frac{2 a^3}{1 + a^2} \check{I}_{22}[\hat{k}_{21,i}](\delta_1^n \sigma_4(x,y),0) \nonumber\\
		&\qquad- \delta_1^n \delta_2^n \frac{2a^2}{\sqrt{\varepsilon_1}(1+a^2)} (\check{\Gamma}_2 \circ g)[k_{21,i}] (\delta_1^n \sigma_4(x,y)) \bigg] \nonumber\\
		&\quad+ \hat{I}_{21}[\hat{k}_{21,i}](x,y) \label{eq:k2_iter_def}
	\end{align}
	We seek the existence fixed point of the iteration (in the complete space $C(\mathcal{T})$), that is, $\lim_{i \rightarrow \infty} \hat{k}_{21,i} = \hat{k}_{21}$. For the limit to converge, we can seek for \emph{uniform} bounds over $\mathcal{T}$, which will imply uniform convergence. For simplicity of analysis, we study the residuals of the iteration, $\Delta \hat{k}_{21,i} := \hat{k}_{21,i+1} - \hat{k}_{21,i}$. Noting the linearitiy of the integral operators, and evaluating the infinite geometric sum, one can eventually find the following iterative bound:
	\begin{align}
		&|\Delta \hat{k}_{21,i+1}(x,y)| \leq \nonumber\\
		&\quad \frac{1}{1- \delta_1 \delta_2} \bigg[|\hat{I}_{21}[|\Delta \hat{k}_{21,i}|](\delta_3^n \sigma_4(x,y), \delta_3^n \sigma_4(x,y))| \nonumber\\
		&\qquad\qquad\qquad+ \delta_2 |\check{I}_{21}[|\Delta \hat{k}_{21,i}|](\delta_1^n \sigma_4(x,y),0)| \nonumber\\
		&\qquad\qquad\qquad+ \frac{2 a^3}{1 + a^2} |\check{I}_{22}[|\Delta \hat{k}_{21,i}|](\delta_1^n \sigma_4(x,y),0)| \nonumber\\
		&\qquad\qquad\qquad+ \frac{2a^2}{\sqrt{\varepsilon_1}(1+a^2)} \int_0^{\delta_1^n \sigma_4(x,y)} \frac{\sqrt{\varepsilon_1} - \sqrt{\varepsilon_2}}{\varepsilon_2} \nonumber\\
		&\qquad\qquad\qquad\quad \times \norm{\Gamma_1}_{L^\infty} \norm{S}_{L^\infty} |\Delta \hat{k}_{21,i}(\xi,\xi)| d\xi \bigg] \nonumber\\
		&\quad+ |\hat{I}_{21}[|\Delta \hat{k}_{21,i}|](x,y)| \label{eq:delta2_bound}
	\end{align}
	where now we can establish the postulated fixed point as
	\begin{align}
		\hat{k}_{21} &= \lim_{i \rightarrow \infty} \hat{k}_{21,i} = \hat{k}_{21,0} + \sum_{i=0}^{\infty} \Delta \hat{k}_{21,i} \label{eq:delta_k2}
	\end{align}
	The equivalent condition for convergence of the residuals is showing that the sum in \eqref{eq:delta_k2} converges uniformly, which we will show via the Weierstrass M-test. We select $\hat{k}_{21,0} = 0$, and from \eqref{eq:k2_iter_def}, we can find
	\begin{align}
		\Delta \hat{k}_{21,0} \leq \frac{1}{1- \delta_1 \delta_2} \frac{2 a^2}{\sqrt{\varepsilon_1}(1+a^2)}|\Gamma_0|_2 \norm{S}_{L^\infty} =: \Psi_0
	\end{align}
	Noting this, one may substitute $\Psi_0$ into \eqref{eq:delta2_bound} to compute successive bounds on $\Delta \hat{k}_{21,i}$, which can be encapsulated into the following bound:
	\begin{align}
		|\Delta k_{21,i}(x,y)| \leq 3\Psi_0 \frac{1}{n!} \Psi^n x^n \label{eq:deltak2_bound_n}
	\end{align}
	where
	\begin{align}
		\Psi &= \max \bigg\{ \left(\frac{1}{1- \delta_1 \delta_2}\right)\left( 1 + \delta_2 + \frac{2 a^3}{1+a^2} \right) \nonumber\\
		&\qquad\qquad\times \left(\frac{\sqrt{\varepsilon_1} - \sqrt{\varepsilon_2}}{\sqrt{\varepsilon_2}}\norm{k_{11}}_{L^\infty} \right), \nonumber\\
		&\qquad\qquad \left(\frac{1}{1- \delta_1 \delta_2}\right) \left( \frac{2}{\varepsilon_2} \norm{\Gamma_1}_{L^\infty} \norm{S}_{L^\infty} \right), \nonumber\\
		&\qquad\qquad \frac{\sqrt{\varepsilon_1}-\sqrt{\varepsilon_2}}{\sqrt{\varepsilon_2}} \norm{k_{11}}_{L^\infty} \bigg\}
	\end{align}
	Then it is quite clear that
	\begin{align}
		\sum_{i = 0}^\infty |\Delta k_{21,i}(x,y)| \leq 3 \Psi_0 \exp(\Psi x) < \infty
	\end{align}
	since $x \in [0,1]$. By the Weierstrass M-test we can conclude the uniform (and absolute) convergence of \eqref{eq:delta_k2} in $C(\mathcal{T})$. This proves the existence of $\hat{k}_{21} \in C(\mathcal{T})$ via the Schauder fixed point theorem.

	As mentioned previously, the existence of $\hat{k}_{21} \in C(\mathcal{T})$ will imply the existence of $\hat{k}_{22}, \check{k}_{21},\check{k}_{22} \in C(\mathcal{T})$ by a mere straightforward evaluation of \eqref{eq:kh_22_eq},\eqref{eq:kc_21_eq},\eqref{eq:kc_22_eq}.

\end{pf}

\begin{lem}
	The ODE \eqref{eq:kernel_ode} and associated initial conditions admit a unique $C([0,1])$ solution.
	\label{lem:wellposedode}
\end{lem}
\begin{pf}
	The proof directly follows from \eqref{eq:gamma1_soln},\eqref{eq:gamma2_soln} and Lemmas \ref{lem:wellposedk1},\ref{lem:wellposedk2}. It is trivial to see that $\Gamma_1 \in C^\infty([0,1])$ by virtue of the matrix exponential.

	The regularity of $\Gamma_2$ can be recovered noting that $\Gamma_2$ involves a convolution of the operator $g[k_{21}](x)$ with a matrix exponential (seen in \eqref{eq:gamma2_soln}). As the exponential is $C^\infty([0,1])$, it is quite clear that it acts as a mollifier to recover a $C^\infty([0,1])$ solution for $\Gamma_2$.
\end{pf}

\section{Conclusion}
\label{sec:conclusion}
A control design methodology via folding and infinte-dimensional backstepping is detailed in the paper. The result comes as a natural extension to the folding framework to designing bilateral controllers.

Of great interest is an alternative interpretation of the folded system with an ODE. The control of an ODE through two distinct controllers can be seen to play a cooperative ``game'' (though not in an optimization sense), whose objective is to stabilize the ODE through the coupled actuation path. This interpretation naturally raises the question of casting a noncooperative game, where perhaps the two controllers are designed independently of one another. Such a formulation may lead to more robust bilateral implementations, where the failure of one controller does not compromise the stability of the system.

A nautral extension to consider is state estimation. Some work by \cite{camacho2019pdeode} has explored a similar problem, albeit for a single parabolic equation (as opposed to two distinct parabolic input paths). The state estimation analogue to the problem has been considered without the ODE in (...), where two collocated measurements of state and flux are taken at an arbitrary interior point (independent of the ODE coupling/folding point). When only the ODE is measured, however, the designer can really only generate an estimate of the state at the point of coupling, and not necessarily of the flux -- an undersensed system. This problem is of great engineering interest, and is under investigation.

\bibliographystyle{ifacconf}
\bibliography{ref2}

\end{document}